\documentclass[reqno,11pt]{amsart}
 
\usepackage{amsmath,amsthm}
\usepackage{amssymb}
\usepackage{euscript}

\numberwithin{equation}{subsection}  

\newcommand{\sqsp}{\renewcommand{\baselinestretch}{1.09}\tiny\normalsize}

\newtheorem{thm}[subsection]{Theorem}
\newtheorem{lemma}[subsection]{Lemma}

\newtheorem{cor}[subsection]{Corollary}

\newtheorem*{conjA}{Conjecture A}
\newtheorem*{conjB}{Conjecture B}
\newtheorem*{conjC}{Conjecture C}

\newcommand{\bZ}{\mathbf{Z}}
\newcommand{\bH}{\mathbf{H}}

\DeclareMathOperator{\genus}{Genus}
\DeclareMathOperator{\lcm}{lcm}
\DeclareMathOperator{\Id}{Id}

\begin{document}
\title{On $\lambda$-rings and topological realization}
\author{Donald Yau}

\begin{abstract}
It is shown that most possibly truncated power series rings admit uncountably many filtered $\lambda$-ring structures.  The question of how many of these filtered $\lambda$-ring structures are topologically realizable by the $K$-theory of torsionfree spaces is also considered for truncated polynomial rings. 
\end{abstract}
\email{dyau@math.uiuc.edu}
\address{Department of Mathematics, University of Illinois at Urbana-Champaign, 1409 W. Green Street, Urbana, IL 61801}
\date{\today}
\maketitle
\sqsp


\section{Introduction}

A $\lambda$-ring is, roughly speaking, a commutative ring $R$ with unit together with operations $\lambda^i$, $i \geq 0$, on it that act like the exterior power operations.  It is widely used in Algebraic Topology, Algebra, and Representation Theory.  For example, the complex representation ring $R(G)$ of a group $G$ is a $\lambda$-ring, where $\lambda^i$ is induced by the map that sends a representation to its $i$th exterior power.  Another example of a $\lambda$-ring is the complex $K$-theory of a topological space $X$.  Here, $\lambda^i$ arises from the map that sends a complex vector bundle $\eta$ over $X$ to the $i$th exterior power of $\eta$.  In the algebra side, the universal Witt ring $\mathbf{W}(R)$ of a commutative ring $R$ is a $\lambda$-ring.

The purpose of this paper is to consider the following two inter-related questions:
   \begin{itemize}
   \item Classify the $\lambda$-ring structures over power series and truncated polynomial rings.
   \item Which ones and how many of these $\lambda$-rings are realizable as (i.e.\ isomorphic to) the $K$-theory of a topological space?
   \end{itemize}
The first question is purely algebraic, with no topology involved.  One can think of the second question as a $K$-theoretic analogue of the classical Steenrod question, which asks for a classification of polynomial rings (over the field of $p$ elements and has an action by the mod $p$ Steenrod algebra) that can be realized as the singular mod $p$ cohomology of a topological space.

In addition to being a $\lambda$-ring, the $K$-theory of a space is filtered, making $K(X)$ a filtered $\lambda$-ring.  Precisely, by a \emph{filtered} $\lambda$-\emph{ring} we mean a filtered ring $(R, \lbrace R=I^0 \supset I^1 \supset \cdots \rbrace)$ in which $R$ is a $\lambda$-ring and the filtration ideals $I^n$ are all closed under the $\lambda$-operations $\lambda^i$ $(i > 0)$.  It is, therefore, more natural for us to consider filtered $\lambda$-ring structures over filtered rings.   Moreover, we will restrict to \emph{torsionfree} spaces, i.e.\ spaces whose integral cohomology is $\bZ$-torsionfree.  The reason for this is that one has more control over the $K$-theory of such spaces than non-torsionfree spaces.

A discussion of our main results follow.  Proofs are mostly given in later sections.  Our first result shows that there is a huge diversity of filtered $\lambda$-ring structures over most power series and truncated polynomial rings.

\medskip
\begin{thm}
\label{thm:existence}
Let $x_1, \ldots, x_n$ be algebraically independent variables with the same (arbitrary but fixed) filtration $d > 0$, and let $r_1, \ldots, r_n$ be integers $\geq 2$, possibly $\infty$.  Then the possibly truncated power series filtered ring $\bZ \lbrack \lbrack x_1, \ldots, x_n \rbrack \rbrack/(x_1^{r_1}, \ldots, x_n^{r_n})$ admits uncountably many isomorphism classes of filtered $\lambda$-ring structures.
\end{thm}

Here $x_i^\infty$ is by definition equal to $0$.  In particular, this Theorem covers both finitely generated power series rings and truncated polynomial rings.  The case $r_i = \infty$ $(1 \leq i \leq n)$ is proved in \cite{yau2}.  In this case, there are uncountably many isomorphism classes that are topologically realizable, namely, by the spaces in the localization genus of $(BS^3)^{\times n}$.  The remaining cases are proved by directly constructing uncountably many filtered $\lambda$-rings.  In general, there is no complete classification of all of the isomorphism classes of filtered $\lambda$-ring structures.  However, such a classification can be obtained for small truncated polynomial rings, in which case we can also give some answers to the second question above.  This will be considered below after a brief discussion of Adams operations.


\subsection*{Adams operations}
\label{subsec:adams op}

The results below are all described in terms of Adams operations.  We will use a result of Wilkerson \cite{wil} on recovering the $\lambda$-ring structure from the Adams operations.  More precisely, Wilkerson's Theorem says that if $R$ is a $\bZ$-torsionfree ring which comes equipped with ring endomorphisms $\psi^k$ $(k \geq 1)$ satisfying the conditions, (1) $\psi^1 = \Id$ and $\psi^k\psi^l = \psi^{kl} = \psi^l \psi^k$ and (2) $\psi^p(r) \equiv r^p \pmod{pR}$ for each prime $p$ and $r \in R$, then $R$ admits a unique $\lambda$-ring structure with the $\psi^k$ as the Adams operations.  The obvious filtered analogue of Wilkerson's Theorem is also true for the filtered rings considered in Theorem \ref{thm:existence}.  Therefore, for these filtered rings, in order to describe a filtered $\lambda$-ring structure, it suffices to describe the Adams operations.

When there are more than one $\lambda$-rings in sight, we will sometimes write $\psi^n_R$ to denote the Adams operation $\psi^n$ in $R$.


\subsection*{Truncated polynomial rings}
\label{subsec:truncated}

Consider the filtered truncated polynomial rings $(\bZ \lbrack x \rbrack/(x^n), \vert x \vert = d)$ and $(\bZ \lbrack x \rbrack/(x^n), \vert x \vert = d^\prime)$, with $x$ in filtration exactly $d > 0$ and $d^\prime > 0$, respectively.  Let $\Lambda_d$ denote the set of isomorphism classes of filtered $\lambda$-ring structures over $(\bZ \lbrack x \rbrack/(x^n), \vert x \vert = d)$.  Define $\Lambda_{d^\prime}$ similarly.  Then it is obvious that $\Lambda_d$ and $\Lambda_{d^\prime}$ are in one-to-one correspondence, and there is no reason to distinguish between them.  Indeed, for $R \in \Lambda_d$ one can associate to it $R^\prime \in \Lambda_{d^\prime}$ such that $\psi^k_R(x) = \psi^k_{R^\prime}(x)$ as polynomials for all $k$, and this construction gives the desired bijection between $\Lambda_d$ and $\Lambda_{d^\prime}$.  Therefore, using the above bijections, we will identify the sets $\Lambda_d$ for $d = 1, 2, \ldots$, and write $\Lambda(\bZ \lbrack x \rbrack/(x^n))$ for the identified set.  Each isomorphism class of filtered $\lambda$-ring structures on $\bZ \lbrack x \rbrack/(x^n)$ with $\vert x \vert \in \lbrace 1, 2, \ldots \rbrace$ is considered an element in $\Lambda(\bZ \lbrack x \rbrack/(x^n))$.

When there is no danger of confusion, we will sometimes not distinguish between a filtered $\lambda$-ring structure and its isomorphism class.

We start with the simplest case $n = 2$, i.e.\ the dual number ring $\bZ \lbrack x \rbrack/(x^2)$.


\medskip
\begin{thm}[= Corollary 4.1.2 in \cite{yau1}]
\label{thm:n=2}
There is a bijection between $\Lambda(\bZ \lbrack x \rbrack/(x^2))$ and the set of sequences $(b_p)$ indexed by the primes in which the component $b_p$ is divisible by $p$.  The filtered $\lambda$-ring structure corresponding to $(b_p)$ has Adams operations $\psi^p(x) = b_p x$.

Moreover, such a filtered $\lambda$-ring is isomorphic to the $K$-theory of a torsionfree space if and only if there exists an integer $k \geq 1$ such that $b_p = p^k$ for all $p$.
\end{thm}

Therefore, in this case, exactly countably infinitely many (among the uncountably many) isomorphism classes are topologically realizable by torsionfree spaces.  Indeed, the $K$-theory of the even-dimensional sphere $S^{2k}$ realizes the filtered $\lambda$-ring with $b_p = p^k$ for all $p$.  On the other hand, if $X$ is a torsionfree space with $K(X) = \bZ \lbrack x \rbrack/(x^2)$, then $\psi^p(x) = p^kx$ for all $p$ if $x$ lies in filtration exactly $2k$ \cite[Corollary 5.2]{adams}.


In what follows we will use the notation $\theta_p(n)$ to denote the largest integer for which $p^{\theta_p(n)}$ divides $n$, where $p$ is any prime.  By convention we set $\theta_p(0) = -\infty$.

To study the case $n = 3$, we need to consider the following conditions for a sequence $(b_p)$ of integers indexed by the primes:  
   \begin{description}
   \item[(A)] $b_2 \not= 0$, $b_p \equiv 0 \pmod{p}$ for all primes $p$, and $b_p(b_p - 1) \equiv 0 \pmod{2^{\theta_2(b_2)}}$ for all odd primes $p$.
   \item[(B)] Let $(b_p)$ be as in (A).  Consider a prime $p > 2$ (if any) for which $b_p \not= 0$ and $\theta_p(b_p) = \min\lbrace \theta_p(b_q(b_q - 1)) \colon \, b_q \not= 0 \rbrace$. 
   \end{description}

We will call them conditions (A) and (B), respectively.  Notice that in (B), $\theta_p(b_p) = \theta_p(b_p(b_p - 1))$, since $p$ does not divide $(b_p - 1)$.  Moreover,  there are at most finitely many such primes, since each such $p$ divides $b_2(b_2 - 1)$, which is non-zero.

The following result gives a complete classification for $\Lambda(\bZ \lbrack x \rbrack/(x^3))$, the set of isomorphism classes of filtered $\lambda$-ring structures over the filtered truncated polynomial ring $\bZ \lbrack x \rbrack/(x^3)$.  As before we will describe a $\lambda$-ring in terms of its Adams operations.

\medskip
\begin{thm}
\label{thm:n=3}
Let $R$ be a filtered $\lambda$-ring structure on $\bZ \lbrack x \rbrack/(x^3)$.  Then $R$ is isomorphic to one of the following filtered $\lambda$-rings:
   \begin{enumerate}
   \item $S((c_p)) = \lbrace \psi^p(x) = c_p x^2 \colon \, p \text{ prime} \rbrace$ with $c_2 \equiv 1 \pmod{2}$ and $c_p \equiv 0 \pmod{p}$ for $p > 2$.   Moreover, any such sequence $(c_p)$ gives rise to an element of $\Lambda(\bZ \lbrack x \rbrack/(x^3))$, and two such filtered $\lambda$-rings, $S((c_p))$ and $S((c^\prime_p))$, are isomorphic if and only if $(c_p) = \pm (c^\prime_p)$.
\medskip
   \item  $S((b_p), k) = \lbrace \psi^p(x) = b_px + c_px^2 \colon \,  p \text{ prime} \rbrace$ with $(b_p)$ satisfying condition (A) above and the $c_p$ having the following form.  Let $p_1, \ldots, p_n$ be the list of all odd primes satisfying condition (B) above.  Then there exists an odd integer $k$ such that, writing $G$ for $\gcd(b_p(b_p - 1) \colon \text{all primes } p)$, 
     \begin{itemize}
     \item $1 \leq k \leq G/2$, 
     \item $k \equiv 0 \pmod{p_1\cdots p_n}$, and
     \item $c_p = kb_p(b_p - 1)/G$ for all primes $p$.
     \end{itemize}
Moreover, any such pair $((b_p), k)$ gives rise to an element of $\Lambda(\bZ \lbrack x \rbrack/(x^3))$.  Two such filtered $\lambda$-rings, $S((b_p), k)$ and $S((b_p^\prime), k^\prime)$, are isomorphic if and only if $b_p = b_p^\prime$ for all primes $p$ and $k = k^\prime$.  No $S((c_p))$ is isomorphic to any $S((b_p), k)$.
   \end{enumerate}
\end{thm}

Notice that in this Theorem, there are uncountably many isomorphism classes in each of cases (1) and (2).  Also note that if there is no prime $p$ satisfying condition (B), then $p_1 \cdots p_n$ is the empty product (i.e.\ $1$), and $k \equiv 0 \pmod{p_1 \cdots p_n}$ is an empty condition.

In case (2), if $(b_p)$ satisfies condition (A), then it follows that there are exactly 
   \begin{equation}
   \label{eq:n=3 bound}
   \left\lceil \frac{G}{4p_1 \cdots p_n} \right\rceil 
   ~=~ \left\lceil \frac{\gcd(b_p(b_p - 1)\colon \text{all } p)}{4p_1 \cdots p_n} \right\rceil 
   \end{equation}
isomorphism classes of filtered $\lambda$-ring structures over $\bZ \lbrack x \rbrack/(x^3)$ with the property that $\psi^p(x) \equiv b_p x \pmod{x^2}$ for all primes $p$.  Here $\lceil s \rceil$ denotes the smallest integer that is greater than or equal to $s$.  Applying formula \eqref{eq:n=3 bound} to the cases $(b_p) = (p^r)$, where $r \in \lbrace 1, 2, 4 \rbrace$, we see that there are a total of $64$ elements in $\Lambda(\bZ \lbrack x \rbrack/(x^3))$ satisfying $\psi^p(x) = p^r x \pmod{x^2}$ for all primes $p$.  In fact, there is a unique such element when $r = 1$ $(k = 1)$, three such elements when $r = 2$ $(k \in \lbrace 1, 3, 5 \rbrace)$, and sixty such elements when $r = 4$ $(k \in \lbrace 1, 3, \ldots , 119\rbrace)$.  This simple consequence of Theorem \ref{thm:n=3} leads to the following upper bound for the number of isomorphism classes of filtered $\lambda$-ring structures over $\bZ \lbrack x \rbrack/(x^3)$ that are topologically realizable by torsionfree spaces.

\pagebreak
\begin{cor}
\label{cor:n=3}
Let $X$ be a torsionfree space whose $K$-theory filtered ring is $\bZ\lbrack x \rbrack/(x^3)$.  Then, using the notation of Theorem \ref{thm:n=3}, $K(X)$ is isomorphic as a filtered $\lambda$-ring to $S((p^r), k)$ for some $r \in \lbrace 1, 2, 4 \rbrace$ and some $k$.  In particular, at most $64$ of the uncountably many isomorphism classes of filtered $\lambda$-ring structures on $\bZ \lbrack x \rbrack/(x^3)$ can be topologically realized by torsionfree spaces.
\end{cor}

Indeed, if $X$ is a torsionfree space whose $K$-theory filtered ring is $\bZ \lbrack x \rbrack/(x^3)$, then by Adams's result on Hopf invariant $1$ \cite{adams1}, the generator $x$ must have filtration exactly $2$, $4$, or $8$.  When the filtration of $x$ is equal to $2r$, one has $\psi^p(x) \equiv p^r x \pmod{x^2}$ for all primes $p$.  Therefore, by Theorem \ref{thm:n=3} $K(X)$ must be isomorphic to $S((p^r), k)$ for some $r \in \lbrace 1, 2, 4 \rbrace$ and some $k$.  The discussion preceding this Corollary then implies that there are exactly $64$ such isomorphism classes of filtered $\lambda$-rings.

It should be remarked that at least $3$ of these $64$ isomorphism classes are actually realized by spaces, namely, the projective $2$-spaces $FP^2$, where $P$ denotes the complex numbers, the quaternions, or the Cayley octonions.  These spaces correspond to $r = 1, 2,$ and $4$, respectively.  Further work remains to be done to determine whether any of the other $61$ filtered $\lambda$-rings in Corollary \ref{cor:n=3} are topologically realizable.

Before moving on to the case $n = 4$, we would like to present another topological application of Theorem \ref{thm:n=3}, which involves the notion of Mislin genus.  Let $X$ be a nilpotent space of finite type (i.e.\ its homotopy groups are all finitely generated and $\pi_n(X)$ is a nilpotent $\pi_1(X)$-module for each $n \geq 2$).  The \emph{Mislin genus} of $X$, denoted $\genus(X)$, is the set of homotopy types of nilpotent finite type spaces $Y$ such that the $p$-localizations of $X$ and $Y$ are homotopy equivalent for all primes $p$.  Let $\bH$ denote the quaternions.  It is known that $\genus(\bH P^\infty)$ is an uncountable set \cite{rector}.  Moreover, these uncountably many homotopically distinct spaces have isomorphic $K$-theory filtered rings \cite{yau2} but pairwise non-isomorphic $K$-theory filtered $\lambda$-rings \cite{notbohm}.  In other words, $K$-theory filtered $\lambda$-ring classifies the Mislin genus of $\bH P^\infty$.  It is also known that $\genus(\bH P^2)$ has exactly $4$ elements \cite{mcgibbon}.  We will now show that the genus of $\bH P^2$ behaves very differently from that of $\bH P^\infty$ as far as $K$-theory filtered $\lambda$-rings are concerned.

\medskip
\begin{cor}
\label{cor2:n=3}
$K$-theory filtered $\lambda$-ring does not classify the Mislin genus of $\bH P^2$.  In other words, there exist homotopically distinct spaces $X$ and $Y$ in $\genus(\bH P^2)$ whose $K$-theory filtered $\lambda$-rings are isomorphic.
\end{cor}

Indeed, an argument similar to the one in \cite{yau2} shows that the $4$ homotopically distinct spaces in the Mislin genus of $\bH P^2$ all have $\bZ \lbrack x \rbrack/(x^3)$ as their $K$-theory filtered ring, with $x$ in filtration exactly $4$.  Therefore, in each one of these $\lambda$-rings, we have $\psi^p(x) \equiv p^2 x \pmod{x^2}$ for all primes $p$.  The Corollary now follows, since by Theorem \ref{thm:n=3}, there are only $3$ isomorphism classes of filtered $\lambda$-ring structures on $\bZ \lbrack x \rbrack/(x^3)$ of the form $S((p^2), k)$ because $k$ must be $1$, $3$, or $5$.  It is still an open question as to whether $K$-theory filtered $\lambda$-ring classifies the genus of $\bH P^n$ for $2 < n < \infty$.


We now move on to the case $n = 4$, i.e.\ the filtered truncated polynomial ring $\bZ \lbrack x \rbrack/(x^4)$.  A complete classification theorem along the lines of Theorem \ref{thm:n=2} and Theorem \ref{thm:n=3} has not yet been achieved for $n \geq 4$.  However, some sort of classification is possible if one imposes certain conditions on the linear coefficients of the Adams operations that usually appear in the $K$-theory of spaces.

\medskip
\begin{thm}
\label{thm:n=4}
Let $R$ be an element of $\Lambda(\bZ \lbrack x \rbrack/(x^4))$.  
   \begin{enumerate}
   \item If $\psi^p_R(x) \equiv px \pmod{x^2}$ for all primes $p$, then $R$ is isomorphic to the  filtered $\lambda$-ring structure with $\psi^p(x) = (1 + x)^p - 1$ for all primes $p$.
\medskip
   \item If $\psi^p_R(x) \equiv p^2x \pmod{x^2}$ for all primes $p$, then $R$ is isomorphic to one of the following $60$ mutually non-isomorphic filtered $\lambda$-ring structures on $\bZ \lbrack x \rbrack/(x^4)$:
     \[
     S(k, d_2) ~=~ \left\lbrace \psi^p(x) = p^2x ~+~ \frac{kp^2(p^2 - 1)}{12}x^2 ~+~ d_px^3 \right\rbrace,
     \]
where $k \in \lbrace 1, 5 \rbrace$, $d_2 \in \lbrace 0, 2, 4, \ldots, 58 \rbrace$, and 
     \begin{equation}
     \label{eq:d_p}
     d_p ~=~ \frac{p^2(p^4 - 1)d_2}{60} 
             ~+~ \frac{k^2p^2(p^2 - 1)(p^2 - 4)}{360}
     \end{equation}
for odd primes $p$.
   \item In general, if $\psi^p_R(x) \equiv b_p x \pmod{x^2}$ with $b_p \not= 0$ for all primes $p$, then there are only finitely many isomorphism classes of filtered $\lambda$-ring structures $S$ over $\bZ \lbrack x \rbrack/(x^4)$ such that $\psi^p_S(x) \equiv b_p x \pmod{x^2}$ for all primes $p$.
   \item If (in the notation of the previous statement) $b_2 = 0$, then $R$ is  of the form $S((c_p), (d_p)) = \lbrace \psi^p(x) = c_p x^2 + d_p x^3 \rbrace$.  Any such collection of polynomials gives rise to a filtered $\lambda$-ring structure, provided that $\psi^p(x) \equiv x^p \pmod{p}$ for all primes $p$.  Two such filtered $\lambda$-rings, $S((c_p), (d_p))$ and $S((\bar{c}_p), (\bar{d}_p))$, are isomorphic if and only if (i) $(c_p) = \pm (\bar{c}_p)$ and (ii) there exists an integer $\alpha$ such that $\bar{d}_p = d_p + 2c_p \alpha$ for all primes $p$.
   \end{enumerate}
\end{thm}

The first two statements of this Theorem immediately leads to the following upper bound for the number of isomorphism classes of filtered $\lambda$-ring structures on $\bZ \lbrack x \rbrack/(x^4)$ that can be topologically realized by torsionfree spaces.

\medskip
\begin{cor}
\label{cor:n=4}
Let $X$ be a torsionfree space whose $K$-theory filtered ring is $\bZ \lbrack x \rbrack/(x^4)$.  Then, as a filtered $\lambda$-ring, $K(X)$ is isomorphic to one of the filtered $\lambda$-rings described in parts (1) and (2) in Theorem \ref{thm:n=4}.  In particular, at most $61$ of the uncountably many isomorphism classes of filtered $\lambda$-ring structures on $\bZ \lbrack x \rbrack/(x^4)$ can be realized as the $K$-theory of a torsionfree space.
\end{cor}

Indeed, if $X$ is such a space, then the filtration of $x$ must be exactly $2$ or $4$ \cite[Corollary 4L.10]{hatcher}, and therefore the linear coefficient of $\psi^p(x)$ for any prime $p$ must be $p$ (if the filtration of $x$ is $2$) or $p^2$ (if the filtration of $x$ is $4$).  So the result follows immediately from Theorem \ref{thm:n=4}.

It should be noted that at least $2$ of the $61$ isomorphism classes in the first two statements of Theorem \ref{thm:n=4} are topologically realizable, namely, by the projective $3$-spaces $FP^3$, where $F$ is either the complex numbers or the quaternions.  The $K$-theory of the former space is case (1) in Theorem \ref{thm:n=4}, while the latter space has $S(1,0)$ as its $K$-theory.  It is still an open question as to whether any of the other $59$ isomorphism classes are topologically realizable.


What happens when $n > 4$, as far as the two questions stated in the beginning of this Introduction are concerned, are only partially understood.  We will discuss several conjectures and some partial results in this general setting.

Concerning the problem of topological realizations, we believe that the finiteness phenomenon in Corollary \ref{cor:n=3} and Corollary \ref{cor:n=4} should not be isolated examples.


\medskip
\begin{conjA}
\label{conjA}
Let $n$ be any integer $\geq 3$.  Then, among the uncountably many isomorphism classes of filtered $\lambda$-ring structures over the filtered truncated polynomial ring $\bZ \lbrack x \rbrack/(x^n)$, only finitely many of them can be realized as the $K$-theory of torsionfree spaces.
\end{conjA}

As mentioned above, the cases $n = 3$ and $4$ are known to be true.  Just like the way Corollary \ref{cor:n=3} and Corollary \ref{cor:n=4} are proved, one way to approach this Conjecture is to consider $\lambda$-rings with given linear coefficients in its Adams operations.  More precisely, we offer the following conjecture.


\medskip
\begin{conjB}
\label{conjB}
Let $n$ be any integer $\geq 3$ and let $\lbrace b_p \in p \bZ \colon p \text{ prime} \rbrace$ be non-zero integers.  Then there exist only finitely many isomorphism classes of filtered $\lambda$-ring structures on $\bZ \lbrack x \rbrack/(x^n)$ with the property that $\psi^p(x) \equiv b_p x \pmod{x^2}$ for all primes $p$.
\end{conjB}

In fact, Conjecture A for $n \geq 4$ would follow from the cases, $(b_p = p)$ and $(b_p = p^2)$, of Conjecture B.  The $n = 3$ case of Conjecture B is contained in Theorem \ref{thm:n=3}.  Moreover, \eqref{eq:n=3 bound} gives the exact number of isomorphism classes in terms of the $b_p$.  One plausible way to prove this Conjecture is to consider extensions of $\lambda$-ring structures one degree at a time.


\medskip
\begin{conjC}
\label{conjC}
Let $n$ be any integer $\geq 3$ and let $R$ and $S$ be isomorphic filtered $\lambda$-ring structures on $\bZ \lbrack x \rbrack/(x^n)$.  Then there exists an isomorphism $\sigma \colon R \xrightarrow{\cong} S$ with the following property: If $\tilde{R}$ is a filtered $\lambda$-ring structure on $\bZ \lbrack x \rbrack/(x^{n+1})$ such that $\psi^p_{\tilde{R}}(x) \equiv \psi^p_R(x) \pmod{x^n}$ for all primes $p$, then 
   \begin{equation}
   \label{eq:conjC}
   \tilde{S} ~\buildrel \text{def} \over =~
   \Bigl\lbrace \psi^p_{\tilde{S}}(x) = (\sigma^{-1} \circ \psi^p_{\tilde{R}} \circ \sigma)(x) \in \bZ \lbrack x \rbrack/(x^{n+1}) \colon p \text{ primes} \Bigr\rbrace
   \end{equation}
is also a filtered $\lambda$-ring structure on $\bZ \lbrack x \rbrack/(x^{n+1})$.
\end{conjC}

Here $\sigma(x)$ is considered a polynomial in both $\bZ \lbrack x \rbrack/(x^{n})$ and $\bZ \lbrack x \rbrack/(x^{n+1})$, and $\sigma^{-1}(x)$ is the (compositional) inverse of $\sigma(x)$ in $\bZ \lbrack x \rbrack/(x^{n+1})$.  In the definition of $\psi^p_{\tilde{S}}(x)$, the symbol $\circ$ means composition of polynomials.

Observe that for any isomorphism $\sigma \colon R \to S$, one has that $(\psi^p_{\tilde{S}} \circ \psi^q_{\tilde{S}})(x) = (\psi^q_{\tilde{S}} \circ \psi^p_{\tilde{S}})(x)$ for all primes $p$ and $q$.  Thus, to prove Conjecture C, one only needs to show that $\psi^p_{\tilde{S}}(x) \equiv x^p \pmod{p}$ for all primes $p$.  Furthermore, if Conjecture C is true, then $\sigma$ induces an isomorphism $\tilde{R} \cong \tilde{S}$.  Denote by $\Lambda_R$ the set of isomorphism classes of filtered $\lambda$-ring structures $\tilde{R}$ over $\bZ \lbrack x \rbrack/(x^{n+1})$ such that $\psi^p_{\tilde{R}}(x) \equiv \psi^p_R(x) \pmod{x^n}$ for all primes $p$, and define $\Lambda_S$ similarly with $S$ replacing $R$.  Since $\psi^p_{\tilde{S}}(x) \equiv \psi^p_S(x) \pmod{x^n}$ for all primes $p$ and any isomorphism $\tau \colon R \to S$, it follows that Conjecture C implies that $\sigma$ induces an embedding $\Lambda_R \hookrightarrow \Lambda_S$.  In particular, Conjecture B would follow from Conjecture C \emph{and} the following finiteness result.

\medskip
\begin{thm}
\label{thm:n>4}
Let $n$ be any integer $\geq 3$, and let $a_{p,i}$ be integers for $p$ primes and $1 \leq i \leq n - 2$ with $a_{p,1} \not= 0$ for every $p$.  Then the number of isomorphism classes of filtered $\lambda$-ring structures $R$ on $\bZ \lbrack x \rbrack/(x^n)$ satisfying
   \begin{equation}
   \label{eq:extension}
   \psi^p_R(x) ~\equiv~ a_{p,1}x + \cdots + a_{p,n-2}x^{n-2} \pmod{x^{n-1}}
   \end{equation}
for all primes $p$, is at most $\min\bigl\lbrace\vert a_{p,1}^{n-1} - a_{p,1} \vert \colon p \text{ primes}\bigr\rbrace$.
\end{thm}

Notice that if Conjecture C is true for $n \leq N$ for some $N$, then, using Theorem \ref{thm:n>4}, one infers that Conjecture B is true for $n \leq N + 1$, which in turn implies Conjecture A for $n \leq N + 1$.  We summarize this in the following diagram.
   \begin{multline*}
   \bigl\lbrace (\text{Conjecture C})_{n \leq N} + \text{Theorem }\ref{thm:n>4} \bigr\rbrace
   ~\Rightarrow~ (\text{Conjecture B})_{n \leq N + 1} \\
   ~\Rightarrow~ (\text{Conjecture A})_{n \leq N + 1}
   \end{multline*}
In view of these implications, even partial results about Conjecture C would be of interest.

It was mentioned above that in order to prove Conjecture C, one only needs to prove the congruence identity, $\psi^p_{\tilde{S}}(x) \equiv x^p \pmod{p}$, for all primes $p$.  In fact, this only needs to be proved for $p < n$, since the following result takes care of the rest.

\medskip
\begin{thm}
\label{thm:p>n}
Let the assumptions and notations be the same as in the statement of Conjecture C.  If $\sigma \colon R \xrightarrow{\cong} S$ is any filtered $\lambda$-ring isomorphism, then $\psi^p_{\tilde{S}}(x) \equiv x^p \pmod{p}$ for all primes $p \geq n$.
\end{thm}

Using this result, one can show that Conjecture C is true for some small values of $n$.  More precisely, we have the following consequence of Theorem \ref{thm:p>n}.

\medskip
\begin{cor}
\label{cor:p>n}
Conjecture C is true for $n = 3, 4$, and $5$.  Therefore, Conjecture B and Conjecture A are true for $n = 3, 4, 5$, and $6$.
\end{cor}

This Corollary will be proved by directly verifying the congruence identity about $\psi^p_{\tilde{S}}(x)$ for $p < n$.  The arguments for the three cases are essentially the same, and it does \emph{not} seem to go through for $n = 6$ (see the discussion after the proof of this Corollary).  A more sophisticated argument seems to be needed to prove Conjecture C in its full generality.


\subsection*{Organization}

The rest of this paper is organized as follows.  The following section gives a brief account of the basics of $\lambda$-rings and Adams operations, ending with the proof of Theorem \ref{thm:existence}.  Proofs of Theorems \ref{thm:n=3} and \ref{thm:n=4} are in the two sections after the following section.  The results concerning the three Conjectures, namely, Theorems \ref{thm:n>4} and \ref{thm:p>n} and Corollary \ref{cor:p>n}, are proved in the final section.


\section{Basics of $\lambda$-rings}
\label{sec:basics}

The reader may refer to the references \cite{at,knutson} for more in-depth discussion of basic properties of $\lambda$-rings.  We should point out that what we call a $\lambda$-ring here is referred to as a ``special'' $\lambda$-ring in \cite{at}.  All rings considered in this paper are assumed to be commutative, associative, and have a multiplicative unit.


\subsection{$\lambda$-rings}
\label{subsec:lambda-rings}

By a $\lambda$-\emph{ring}, we mean a commutative ring $R$ that is equipped with functions
   \[
   \lambda^i \colon R ~\to~ R \quad (i \geq 0),
   \]
called $\lambda$-operations.  These operations are required to satisfy the following conditions.  For any integers $i, j \geq 0$ and elements $r$ and $s$ in $R$, one has:
   \begin{itemize}
   \item $\lambda^0(r) = 1$.
   \item $\lambda^1(r) = r$.
   \item $\lambda^i(1) = 0$ for $i > 1$.
   \item $\lambda^i(r + s) = \sum_{k = 0}^i\, \lambda^k(r)\lambda^{i-k}(s)$.
   \item $\lambda^i(rs) = P_i(\lambda^1(r), \ldots , \lambda^i(r); \lambda^1(s), \ldots, \lambda^i(s))$.
   \item $\lambda^i(\lambda^j(r)) = P_{i,j}(\lambda^1(r), \ldots , \lambda^{ij}(r))$.
   \end{itemize}
The $P_i$ and $P_{i,j}$ are certain universal polynomials with integer coefficients, and they are defined using the elementary symmetric polynomials as follows.  Given the variables $\xi_1, \ldots , \xi_i$ and $\eta_1, \ldots , \eta_i$, let $s_1, \ldots , s_i$ and $\sigma_1, \ldots , \sigma_i$, respectively, be the elementary symmetric functions of the $\xi$'s and the $\eta$'s.  Then the polynomial $P_i$ is defined by requiring that the expression 
$P_i(s_1, \ldots , s_i; \sigma_1, \ldots , \sigma_i)$ 
be the coefficient of $t^i$ in the finite product
   \[
   \prod_{m,n=1}^i\, (1 + \xi_m \eta_n t).
   \]
The polynomial $P_{i,j}$ is defined by requiring that the expression 
$P_{i,j}(s_1, \ldots , s_{ij})$ 
be the coefficient of $t^i$ in the finite product
   \[
   \prod_{l_1 < \cdots < l_j} \, (1 + \xi_{l_1} \cdots \xi_{l_j} t).
   \]
A $\lambda$-ring map is a ring map which commutes with all the $\lambda$-operations.

By a \emph{filtered ring}, we mean a commutative ring $R$ together with a decreasing sequence of ideals
   \[ 
   R ~=~ I^0 \supseteq I^1 \supseteq I^2 \supseteq \cdots.
   \]
A \emph{filtered ring map} $f \colon R \to S$ is a ring map that preserves the filtrations, i.e.\ $f(I^n_R) \subseteq I^n_S$ for all $n$.

A \emph{filtered} $\lambda$-\emph{ring} is a $\lambda$-ring $R$ which is also a filtered ring in which each ideal $I^n$ is closed under $\lambda^i$ for $i \geq 1$.  Suppose that $R$ and $S$ are two filtered $\lambda$-rings.  Then a filtered $\lambda$-ring map $f \colon R \to S$ is a $\lambda$-ring map that also preserves the filtration ideals.


\subsection{Adams operations}
\label{subsec:adams}

There are some very useful operations inside a $\lambda$-ring $R$, the so-called Adams operations
   \[
   \psi^n \colon R \to R \quad (n \geq 1).
   \]
They are defined by the Newton formula:
   \[
   \psi^n(r) - \lambda^1(r)\psi^{n-1}(r) + \cdots + (-1)^{n-1}\lambda^{n-1}(r) \psi^1(r) + (-1)^n n\lambda^n(r) = 0.
   \]
Alternatively, one can also define them using the closed formula
   \[
   \psi^k ~=~ Q_k(\lambda^1, \ldots, \lambda^k).
   \]
Here $Q_k$ is the integral polynomial with the property that 
   \[
   Q_k(\sigma_1, \ldots, \sigma_k) 
   ~=~ x_1^k ~+~ \cdots ~+~ x_k^k,
   \]
where the $\sigma_i$ are the elementary symmetric polynomials of the $x$'s.  The Adams operations have the following properties:
\begin{enumerate}
\label{adams}
\item All the $\psi^n$ are $\lambda$-ring maps on $R$, and they preserve the filtration ideals if $R$ is a filtered $\lambda$-ring.
\item $\psi^1 = \Id$ and $\psi^m \psi^n = \psi^{mn} = \psi^n \psi^m$.
\item $\psi^p(r) \equiv r^p$ (mod $pR$) for each prime $p$ and element $r$ in $R$.
\end{enumerate}
A $\lambda$-ring map $f$ is compatible with the Adams operations, in the sense that $f\psi^n = \psi^n f$ for all $n$.

The following simple observation will be used many times later in this paper.  Suppose that $R$ and $S$ are $\lambda$-rings with $S$ $\bZ$-torsionfree and that $f \colon R \to S$ is a ring map satisfying $f\psi^p = \psi^pf$ for all primes $p$.  Then $f$ is a $\lambda$-ring map.  Indeed, it is clear that $f$ is compatible with all $\psi^n$ by (2) above.  The Newton formula and the $\bZ$-torsionfreeness of $S$ then imply that $f$ is compatible with the $\lambda^n$ as well.

As discussed in the Introduction, Wilkerson's Theorem \cite{wil} says that if $R$ is a $\mathbf{Z}$-torsionfree ring equipped with ring endomorphisms $\psi^n$ $(n \geq 1)$ satisfying conditions (2) and (3) above, then there exists a unique $\lambda$-ring structure on $R$ whose Adams operations are exactly the given $\psi^n$.  In particular, over the possibly truncated power series filtered ring $\bZ \lbrack \lbrack x_1, \ldots, x_n \rbrack \rbrack/(x_1^{r_1}, \ldots, x_n^{r_n})$ as in Theorem \ref{thm:existence}, a filtered $\lambda$-ring structure is specified by power series $\psi^p(x_i)$ without constant terms, $p$ primes and $1 \leq i \leq n$, such that
   \begin{equation}
   \label{eq:adams1}
   \psi^p(\psi^q(x_i)) ~=~ \psi^q(\psi^p(x_i))
   \end{equation}
and
   \begin{equation}
   \label{eq:adams2}
   \psi^p(x_i) ~\equiv~ x_i^p \pmod{p}
   \end{equation}
 for all such $p$ and $i$.


\subsection{Proof of Theorem \ref{thm:existence}}
\label{subsec:existence}

\begin{proof}
Denote by $R$ the possibly truncated power series filtered ring $\bZ \lbrack \lbrack x_1, \ldots, x_n \rbrack \rbrack/(x_1^{r_1}, \ldots, x_n^{r_n})$ as in the statement of Theorem \ref{thm:existence}.

As was mentioned in the Introduction, the case $r_i = \infty$ for all $i$ is proved in \cite{yau2}.  It remains to consider the cases when at least one $r_i$ is finite.

Assume that at least one $r_i$ is finite.  Let $N$ be the maximum of those $r_j$ that are finite.  For each prime $p \geq N$ and each index $j$ for which $r_j < \infty$, choose an arbitrary positive integer $b_{p,j} \in p \bZ$ such that $b_{p,j} \geq r_j$ and $b_{p,j} \not= p$.  There are uncountably many such choices, since $N < \infty$ and there are countably infinitely many choices for $b_{p,j}$ for each $p \geq N$.  Consider the following power series in $R$:
   \begin{equation}
   \label{eq:adams op}
   \psi^p(x_i) ~=~ \begin{cases}
                   (1 + x_i)^{b_{p,i}} - 1 & \text{ if } p \geq N \text{ and } r_i < \infty, \\
                   (1 + x_i)^p - 1         & \text{ otherwise}.\end{cases}
   \end{equation}
Here $p$ runs through the primes and $i = 1, 2, \ldots, n$.  The collection of power series $\lbrace \psi^p(x_i) \colon 1 \leq i \leq n \rbrace$ extends uniquely to a filtered ring endomorphism $\psi^p$ of $R$.

We first claim that these endomorphisms $\psi^p$, $p$ primes, are the Adams operations of a filtered $\lambda$-ring structure $S$ on $R$.  Since $R$ is $\bZ$-torsionfree, by Wilkerson's Theorem \cite{wil}, it suffices to show that 
   \begin{equation}
   \label{eq:commutativity}
   \psi^p \psi^q ~=~ \psi^q \psi^p
   \end{equation}
and that 
   \begin{equation}
   \label{eq:rp}
   \psi^p(r) ~\equiv~ r^p \pmod{pR}
   \end{equation}
for all primes $p$ and $q$ and elements $r \in R$.  Both of these conditions are verified easily using \eqref{eq:adams op}.  Equation \eqref{eq:commutativity} is true because it is true when applied to each $x_i$ and that the $x_i$ are algebra generators of $R$.  Equation \eqref{eq:rp} is true, since it is true for each $r = x_i$.

Now suppose that $\bar{S}$ is another filtered $\lambda$-ring structure on $R$ constructed in the same way with the integers $\lbrace \bar{b}_{p,j}\rbrace$.  (Here again $p$ runs through the primes $\geq N$ and $j$ runs through the indices for which $r_j < \infty$.)  So in $\bar{S}$, $\psi^p(x_i)$ looks just like it is in \eqref{eq:adams op}, except that $b_{p,i}$ is replaced by $\bar{b}_{p,i}$.  Suppose, in addition, that there is a prime $q \geq N$ such that
   \[
   \lbrace b_{q,j} \rbrace
   ~\not= ~ \lbrace \bar{b}_{q,j} \rbrace
   \]
as sets.  We claim that $S$ and $\bar{S}$ are not isomorphic as filtered $\lambda$-rings.

To see this, suppose to the contrary that there exists a filtered $\lambda$-ring isomorphism
   \[
   \sigma \colon S ~\to~ \bar{S}.
   \]
Let $j$ be one of those indices for which $r_j$ is finite.  Then, modulo filtration $2d$, one has
   \[
   \sigma(x_j) ~\equiv~ a_1 x_1 + \cdots + a_n x_n
   \]
for some $a_1, \ldots, a_n \in \bZ$, not all of which are equal to $0$.  If $r_i = \infty$, we set $\bar{b}_{q,i} = q$.  Equating the linear coefficients on both sides of the equation
   \[
   \sigma \psi^q(x_j) ~=~ \psi^q \sigma(x_j),
   \]
one infers that
   \[
   b_{q,j}\cdot \sum\, a_i x_i 
   ~=~ \sum \, a_i \bar{b}_{q,i}x_i.
   \]
If $a_i \not= 0$ (and such an $a_i$ must exist), then
   \[
   b_{q,j} ~=~ \bar{b}_{q,i}
   \]
for some $i$.  In particular, it follows that $\lbrace b_{q,j} \rbrace$ is contained in $\lbrace \bar{b}_{q,j} \rbrace$.  Therefore the two sets are equal by symmetry.  This is a contradiction.

This finishes the proof of Theorem \ref{thm:existence}.
\end{proof}


\section{Proof of Theorem \ref{thm:n=3}}
\label{sec:n=3}

First we need to consider when a collection of polynomials can be the Adams operations of a filtered $\lambda$-ring structure on $\bZ \lbrack x \rbrack/(x^3)$.  We will continue to describe $\lambda$-rings in terms of their Adams operations.

\medskip
\begin{lemma}
\label{lem1:n=3}
A collection of polynomials, $\lbrace \psi^p(x) = b_p x + c_p x^2 \colon p \text{ prime} \rbrace$, in $\bZ \lbrack x \rbrack/(x^3)$ extends to (the Adams operations of) a filtered $\lambda$-ring structure if and only if the following three statements are satisfied:
   \begin{enumerate}
   \item $b_p \equiv 0 \pmod{p}$ for all primes $p$,
   \item $c_2 \equiv 1 \pmod{2}$ and $c_p \equiv 0 \pmod{p}$ for all primes $p > 2$, and
   \item $(b_q^2 - b_q)c_p = (b_p^2 - b_p)c_q$ for all primes $p$ and $q$.
   \end{enumerate}
Now suppose that these conditions are satisfied.  If $b_2 \not= 0$, then $(b_p^2 - b_p) \equiv 0 \pmod{2^{\theta_2(b_2)}}$ for all odd primes $p$.  If $b_2 = 0$, then $b_p = 0$ for all odd primes $p$. 
\end{lemma}

\begin{proof}
The polynomials $\psi^p(x)$ in the statement above extend to a filtered $\lambda$-ring structure on $\bZ \lbrack x \rbrack/(x^3)$ if and only if \eqref{eq:adams1} and \eqref{eq:adams2} are satisfied.  Expanding $\psi^p(\psi^q(x))$ one obtains
   \[
   \begin{split}
   \psi^p(\psi^q(x)) 
   &~=~ b_p(b_qx + c_qx^2) + c_p(b_qx + c_qx^2)^2 \\
   &~=~ (b_pb_q)x + (b_pc_q + b_q^2 c_p)x^2.
   \end{split}
   \]
Using symmetry and equating the coefficients of $x^2$, it follows that \eqref{eq:adams1} in this case is equivalent to
   \begin{equation}
   \label{eq:condition}
   (b_q^2 - b_q)c_p ~=~ (b_p^2 - b_p)c_q.
   \end{equation}
It is clear that \eqref{eq:adams2} is equivalent to conditions (2) and (3) together, since $x^p = 0$ for $p > 2$.

Now assume that statements (1), (2), and (3) are satisfied.  If $b_2 \not= 0$, then the right-hand side of \eqref{eq:condition} when $p = 2$ is congruent to $0$ modulo $2^{\theta_2(b_2)}$, and therefore so is the left-hand side.  The assertion now follows, since $c_2$ is odd.  If $b_2 = 0$, then the right-hand side of \eqref{eq:condition} when $p = 2$ is equal to $0$, and so $b_q^2 = b_q$ for all odd primes $q$, since $c_2 \not= 0$.  But $b_q \not= 1$, so $b_q = 0$.
\end{proof}

Next we need to know when two filtered $\lambda$-ring structures over $\bZ \lbrack x \rbrack/(x^3)$ are isomorphic.

\medskip
\begin{lemma}
\label{lem2:n=3}
Let $S = \lbrace \psi^p(x) = b_p x + c_p x^2 \rbrace$ and $\bar{S} = \lbrace \psi^p(x) = \bar{b}_p x + \bar{c}_p x^2 \rbrace$ (where $p$ runs through the primes) be two filtered $\lambda$-ring structures over $\bZ \lbrack x \rbrack/(x^3)$.  Then $S$ and $\bar{S}$ are isomorphic filtered $\lambda$-rings if and only if the following two conditions are satisfied simultaneously: 
   \begin{enumerate}
   \item $b_p = \bar{b}_p$ for all primes $p$.
   \item 
       \begin{enumerate}
       \item If $b_2 = \bar{b}_2 = 0$, then there exists $u \in \lbrace \pm 1 \rbrace$ such that $c_p = u \bar{c}_p$ for all primes $p$.
       \item If $b_2 = \bar{b}_2 \not= 0$, then there exist $u \in \lbrace \pm 1 \rbrace$ and $a \in \bZ$ such that 
       \[
       ab_2(b_2 - 1) ~=~ c_2 - u\bar{c}_2.
       \]
       \end{enumerate}
   \end{enumerate}
\end{lemma}

\begin{proof}
Suppose that $S$ and $\bar{S}$ are isomorphic, and let $\sigma \colon S \to \bar{S}$ be a filtered $\lambda$-ring isomorphism.  Then 
   \begin{equation}
   \label{eq:sigma}
   \sigma(x) ~=~ u x + ax^2
   \end{equation}
for some $u \in \lbrace \pm 1 \rbrace$ and integer $a$.  Applying the map $\sigma \psi^p$, $p$ any prime, to the generator $x$, one obtains
   \[
   \sigma\psi^p(x) 
   ~=~ ub_p x + (ab_p + c_p)x^2.
   \]
Similarly, one has
   \[
   \psi^p \sigma(x) ~=~ u\bar{b}_p x + (a\bar{b}_p^2 + u\bar{c}_p)x^2.
   \]
Recall from the previous Lemma that $b_2 = 0$ implies $b_p = 0$ for all odd primes $p$.  Therefore, the ``only if'' part now follows by equating the coefficients in the equation $\sigma \psi^p(x) = \psi^p \sigma(x)$.

Conversely, suppose that conditions (1) and (2)(a) in the statement of the Lemma hold.  Then clearly the map $\sigma \colon S \to \bar{S}$ given on the generator by $\sigma(x) = ux$ is the desired isomorphism.

Now suppose that conditions (1) and (2)(b) hold.  The polynomial $\sigma(x) =   u x + ax^2$ extends uniquely to a filtered ring automorphism on $\bZ \lbrack x \rbrack/x^3$.  The calculation in the first paragraph of this proof shows that, if $b_p = 0$ for a certain prime $p$, then $\sigma \psi^p(x) = \psi^p \sigma(x)$.  If $b_p \not= 0$, then \eqref{eq:condition} in the proof of Lemma \ref{lem1:n=3} implies that
   \[
   a ~=~ \dfrac{c_2 - u \bar{c}_2}{b_2^2 - b_2}
     ~=~ \dfrac{c_p - u \bar{c}_p}{b_p^2 - b_p}.
   \]
Therefore, the argument in the first paragraph once again shows that $\sigma \psi^p(x) = \psi^p \sigma(x)$, and $\sigma \colon S \to \bar{S}$ is the desired isomorphism.
\end{proof}

\medskip
\begin{proof}[Proof of Theorem \ref{thm:n=3}]
It follows immediately from Lemma \ref{lem1:n=3} and Lemma \ref{lem2:n=3} that, in the notation of the statement of Theorem \ref{thm:n=3}, the $S((c_p))$ are all filtered $\lambda$-ring structures on $\bZ \lbrack x \rbrack/(x^3)$ and that two of them are isomorphic if and only if the stated conditions hold.  Similar remarks apply to the $S((b_p), k)$.  Also no $S((c_p))$ is isomorphic to an $S((b_p), k)$.  It remains only to show that any filtered $\lambda$-ring structure $R$ on $\bZ \lbrack x \rbrack/(x^3)$ is isomorphic to one of them.

Write $\psi^p(x) = b_p x + c_p x^2$ for the Adams operations in $R$.  If $b_2 = 0$, then so is $b_p$ for each odd prime $p$.  Then $(c_p)$ satisfies condition (2) in Lemma \ref{lem1:n=3}, and $R$ is isomorphic (in fact, equal) to $S((c_p))$.

Suppose that $b_2 \not= 0$.  Then by Lemma \ref{lem1:n=3} the sequence $(b_p)$ satisfies condition (A).  First consider the case when $b_p = 0$ for all odd primes $p$.  In this case, condition (3) in Lemma \ref{lem1:n=3} implies that $c_p = 0$ for all odd primes $p$.  Let $r$ denote the remainder of $c_2$ modulo $b_2(b_2 - 1)$.  Since $r$ is also an odd integer, there exists a unique odd integer $k$ in the range $1 \leq k \leq b_2(b_2 - 1)/2$ that is congruent $\pmod{b_2(b_2 - 1)}$ to either $r$ or $-r$.  By Lemma \ref{lem2:n=3}, $R$ is isomorphic to $S((b_p), k)$.

Finally, consider the case when $b_2 \not= 0$ and there is at least one odd prime $p$ for which $b_p \not= 0$.  If
   \[
   c_2 ~=~ qb_2(b_2 - 1) + r
   \]
for some integers $q$ and $r$ with $0 \leq r < b_2(b_2 - 1)$, then $r$ must be an odd integer, since $b_2$ is even and $c_2$ is odd.  Define 
   \[
   \bar{c}_2 ~=~ \begin{cases} r & \text{ if } 1 \leq r \leq \dfrac{b_2(b_2 - 1)}{2} \\ b_2(b_2 - 1) - r & \text{ if } r > \dfrac{b_2(b_2 - 1)}{2} \end{cases}
   \]
and 
   \[
   \bar{c}_q ~=~ \begin{cases} c_q - qb_q(b_q - 1) & \text{ if }  1 \leq r \leq \dfrac{b_2(b_2 - 1)}{2} \\ (1 + q)b_q(b_q - 1) - c_q & \text{ if }  r > \dfrac{b_2(b_2 - 1)}{2} \end{cases}
   \]
for $q > 2$.  The three conditions (1) - (3) in Lemma \ref{lem1:n=3} are all easily verified for the polynomials $S = \lbrace \psi^p(x) = b_p x + \bar{c}_px^2 \rbrace$, and so $S$ is a filtered $\lambda$-ring structure on $\bZ\lbrack x \rbrack/(x^3)$.  Observe that $1 \leq \bar{c}_2 \leq b_2(b_2 - 1)/2$.  Moreover, by Lemma \ref{lem2:n=3}, $S$ is isomorphic to $R$.  Equation \eqref{eq:condition} implies that, if $p$ is a prime for which $b_p \not= 0$, then 
   \[
   \bar{c}_2 ~\equiv~ 0 \quad \left(\bmod \, \frac{b_2(b_2 - 1)}{\gcd(b_2(b_2 - 1), b_p(b_p - 1))} \right).
   \]
Therefore, we have 
   \[
   \bar{c}_2 ~\equiv~ 0 \quad \left(\bmod \, \lcm\Bigl( \frac{b_2(b_2 - 1)}{\gcd(b_2(b_2 - 1), b_p(b_p - 1))} \colon \text{all primes }p \Bigr) \right).
   \]
Now observe that
   \[
   \lcm\Bigl( \frac{b_2(b_2 - 1)}{\gcd(b_2(b_2 - 1), b_p(b_p - 1))} \colon \text{all primes }p \Bigr)
   ~=~ \frac{b_2(b_2 - 1)}{G}, 
   \]
where 
   \[
   G ~=~ \gcd(b_p(b_p - 1) \colon \text{all primes }p).
   \]
From the last assertion of Lemma \ref{lem1:n=3}, we see that $2^{\theta_2(b_2)}$ is a factor of $G$, and so $b_2(b_2 - 1)/G$ is odd.   If 
   \[
   \bar{c}_2 ~=~ \frac{l b_2(b_2 - 1)}{G}, 
   \]
then it follows that $l$ must be odd with $1 \leq l \leq G/2$ and that 
   \[
   \bar{c}_p ~=~ \frac{l b_p(b_p - 1)}{G}
   \]
for all primes $p$.   Let $p_1, \ldots, p_n$ be the list of all (if any) odd primes for which condition (B) holds for the sequence $(b_p)$.  For each $p = p_i$, $1 \leq i \leq n$, we have
   \[
   \theta_{p}(\bar{c}_{p}) 
   ~=~ \theta_{p}\left(\frac{l b_p(b_p - 1)}{G}\right)
   ~=~ \theta_{p}(l) 
   ~\geq~ 1,
   \]
and so 
   \[
   l ~\equiv~ 0 \pmod{p}.
   \]
It follows that 
   \[
   l ~\equiv~ 0 \pmod{p_1\cdots p_n}.
   \]
We have shown that $R$ is isomorphic to $S = S((b_p), l)$, as desired.
\end{proof}


\section{Proof of Theorem \ref{thm:n=4}}
\label{sec:n=4}


\subsection{Proof of Theorem \ref{thm:n=4} (1)}
\begin{proof}
First we make the following observation.  Consider a collection of polynomials    \[
   R ~=~ \lbrace \psi^p(x) = b_p x + c_p x^2 + d_p x^3 \colon p \text{ primes} \rbrace
   \]
in $\bZ \lbrack x \rbrack/(x^4)$.  Then $R$ extends to (the Adams operations of) a filtered $\lambda$-ring structure if and only if (i) the $b_p$ and $c_p$ satisfy conditions (1) - (3) of Lemma \ref{lem1:n=3}, (ii)
   \begin{equation}
   \label{eq:d}
   (b_q^3 - b_q)d_p
   ~=~ (b_p^3 - b_p)d_q ~+~ 2c_p c_q(b_p - b_q)
   \end{equation}
for all primes $p$ and $q$, and (iii)
   \begin{equation}
   \label{eq:dp3}
   d_p ~\equiv~ \begin{cases} 0 \pmod p & \text{ if } p \not= 3 \\
                              1 \pmod 3 & \text{ if } p = 3. \end{cases}
   \end{equation}
In fact, \eqref{eq:d} and \eqref{eq:dp3} are obtained from \eqref{eq:adams1} and \eqref{eq:adams2}, respectively, by comparing the coefficients of $x^3$.

Now restrict to case (1) of Theorem \ref{thm:n=4}, i.e.\ when $b_p = p$ for all primes $p$.  By condition (3) in Lemma \ref{lem1:n=3}, we have that in $R$,
   \begin{equation}
   \label{eq1:n=4 bp=p}
   \psi^p(x) ~=~ px ~+~ \frac{c_2 p(p-1)}{2}x^2 ~+~ d_p x^3
   \end{equation}
for each prime $p$. Applying \eqref{eq:d} to this present case, we obtain
   \begin{equation}
   \label{eq2:n=4 bp=p}
   d_p ~=~ \frac{1}{6}p(p - 1)\bigl((p + 1)d_2 ~+~ c_2^2(p - 2)\bigr) \in \bZ.
   \end{equation}
Let $S$ denote the filtered $\lambda$-ring structure on $\bZ \lbrack x \rbrack/(x^4)$ with 
   \[
   \psi^p(x) ~=~ (1 + x)^p ~-~ 1
   \]
for all primes $p$.  Define $\sigma(x)$ by
   \[
   \sigma(x) ~=~ \begin{cases}
                 -x + \dfrac{c_2 + 1}{2}x^2 - \dfrac{(c_2 + 1)(c_2 + 2) + d_2}{6}x^3 & \text{ if } c_2 \equiv 0 \pmod 3 \\
                 x + \dfrac{c_2 - 1}{2}x^2 + \dfrac{(c_2 - 1)(c_2 - 2) + d_2}{6}x^3 & \text{ otherwise}. \end{cases}
   \]
Extend $\sigma(x)$ to a filtered ring automorphism on $\bZ \lbrack x \rbrack/(x^4)$.  Using \eqref{eq1:n=4 bp=p} and \eqref{eq2:n=4 bp=p}, it is now straightforward to check that $\sigma$ is compatible with the Adams operations $\psi^p$ in $R$ and $S$ (i.e.\ $\sigma \psi^p_R(x) = \psi^p_S\sigma(x)$), and so it gives an isomorphism $R \cong S$.
\end{proof}


\subsection{Proof of Theorem \ref{thm:n=4} (2)}
\begin{proof}
Using the case $n = 3$ of Corollary \ref{cor:p>n} (i.e.\ Conjecture C for $n = 3$), which will be proved below, one infers that $R$ is isomorphic to some filtered $\lambda$-ring structure $T$ on $\bZ \lbrack x \rbrack/(x^4)$ with the following property: There exists an integer $k \in \lbrace 1, 3, 5 \rbrace$ such that
   \[
   \psi^p_T(x) ~\equiv~ \psi^p_S(x) \pmod{x^3}
   \]
for all primes $p$, where $S = S((p^2), k)$ is as in Theorem \ref{thm:n=3} (2).  In other words, $R$ must be isomorphic to a ``$\lambda$-ring extension'' to $\bZ \lbrack x \rbrack/(x^4)$ of one of the three $S((p^2), k)$.  Therefore, to prove Theorem \ref{thm:n=4} (2), it suffices to prove the following three statements.
   \begin{enumerate}
   \item The sixty $S(k, d_2)$ listed in the statement of Theorem \ref{thm:n=4} are all filtered $\lambda$-ring structures on $\bZ \lbrack x \rbrack/(x^4)$.
   \item Those $60$ filtered $\lambda$-rings are mutually non-isomorphic.
   \item $T$ (and thus $R$) is isomorphic to one of the sixty $S(k, d_2)$.
   \end{enumerate}
Each one of these three statements is dealt with in a Lemma below.

\medskip
\begin{lemma}
\label{lem1:n=4}
Let $k \in \lbrace 1, 5 \rbrace$ be an integer and let $d_2$ be an even integer.  Consider the collection of polynomials
   \[
   S ~=~ \left\lbrace \psi^p(x) = p^2x ~+~ \frac{kp^2(p^2 - 1)}{12}x^2 ~+~ d_px^3 \right\rbrace,
   \]
where $d_p$ for $p > 2$ is determined by $d_2$ via \eqref{eq:d_p}.  Then $S$ is a filtered $\lambda$-ring structures on $\bZ \lbrack x \rbrack/(x^4)$.
\end{lemma}

Thus, this $S$ looks just like one of those $S(k, d_2)$, except that in $S$ we only require $d_2$ to be even.  In particular, this Lemma implies that each $S(k, d_2)$ is indeed a filtered $\lambda$-ring structure on  $\bZ \lbrack x \rbrack/(x^4)$.

\begin{proof}
We only need to check \eqref{eq:d} and \eqref{eq:dp3} in $S$, since it is clear that conditions (1) - (3) of Lemma \ref{lem1:n=3} are satisfied by $b_p = p^2$ and $c_p = kp^2(p^2 - 1)/12$.  The condition \eqref{eq:d} can be verified directly by using \eqref{eq:d_p}, which expresses each $d_p$ in terms of $d_2$.

Now consider \eqref{eq:dp3}.  The case $p = 2$ is true by hypothesis.  If $p = 3$, then it follows from \eqref{eq:d_p} that
   \begin{equation}
   \label{eq:d3}
   d_3 ~=~ 12d_2 + k^2.
   \end{equation}
Since $k$ is either $1$ or $5$, we must have that $d_3 \equiv 1 \pmod 3$.  If $p = 5$, then 
   \[
   d_5 ~=~ 260d_2 + 35k^2 ~\equiv~ 0 \pmod 5.
   \]

Suppose now that $p > 5$.  We will show that both summands of $d_p$ in \eqref{eq:d_p} are divisible by $p$.  First observe that $12 = 2^2\cdot 3$ must divide $(p - 1)(p + 1)$.  If $5$ divides either $(p - 1)$ or $(p + 1)$, then $60$ divides $(p^2 - 1)$.  If $5$ divides neither one of them, then $p ~\equiv~ 2 \text{ or } 3 \pmod 5$.  In either case, we have 
$p^2 + 1 ~\equiv~ 0 \pmod 5$, and therefore $60$ divides $(p^2 - 1)(p^2 + 1) = (p^4 - 1)$.  It follows that the first summand of $d_p$, namely, $p^2(p^4 - 1)d_2/60$, is an integer that is divisible by $p$ (in fact, by $p^2$).  Now we consider the second summand.  Since $p$ is odd, $(p - 1)p(p + 1)$ is divisible by $2^3\cdot 3$, and so $p(p^2 - 1)(p^2 - 4)$ is divisible by $2^3 \cdot 3^2$.  It is also clear that it is divisible by $5$.  Therefore, $2^3\cdot 3^2 \cdot 5 = 360$ divides $p(p^2 - 1)(p^2 - 4)$.  It follows that the second summand of $d_p$ is an integer that is divisible by $p$.  This proves the Lemma.
\end{proof}

\medskip
\begin{lemma}
\label{lem2:n=4}
Let $S$ and $S^\prime$ be two filtered $\lambda$-ring structures on $\bZ \lbrack x \rbrack/(x^4)$ of the form described in Lemma \ref{lem1:n=4}.  Write $k^\prime$ and $d_p^\prime$ for the coefficients of $x^2$ in $\psi^2_{S^\prime}(x)$ and of $x^3$ in $\psi^p_{S^\prime}(x)$, respectively.  Then $S$ and $S^\prime$ are isomorphic filtered $\lambda$-rings if and only if (i) $k = k^\prime$ and (ii) $d_2 \equiv d_2^\prime \pmod{60}$.
\end{lemma}

In particular, it follows from this Lemma that the sixty $S(k, d_2)$ are mutually non-isomorphic and that any $S$ as in Lemma \ref{lem1:n=4} is isomorphic to one of the $S(k, d_2)$.

\begin{proof}
Suppose that $S$ and $S^\prime$ are isomorphic.  Since their reductions modulo $x^3$ are also isomorphic filtered $\lambda$-rings and since $k, k^\prime \in \lbrace 1, 5 \rbrace$, it follows from Theorem \ref{thm:n=3} that $k = k^\prime$.  Let $\sigma \colon S \to S^\prime$ be an isomorphism, and write 
   \[
   \sigma(x) ~=~ ux + \alpha x^2 + \beta x^3,
   \]
so that $u \in \lbrace \pm 1 \rbrace$.  Consider the equality    
   \begin{equation}
   \label{eq:psi2}
   (\psi^2_S \circ \sigma)(x) ~=~ (\sigma \circ \psi^2_{S^\prime})(x).
   \end{equation}
(Recall that $\circ$ means composition of polynomials.)  The coefficients of $x^2$ on both sides of this equation give rise to the equality
   \[
   4 \alpha + k ~=~ uk + 16 \alpha,
   \]
or, equivalently,
   \[
   (1 - u)k ~=~ 12 \alpha.
   \]
Since $3$ does not divide $k$, the only solution is $u = 1$ and $\alpha = 0$.  The coefficients of $x^3$ in \eqref{eq:psi2} then give rise to the equality
   \[
   4 \beta + d_2 ~=~ d_2^\prime + 64 \beta.
   \]
This proves the ``only if'' part of the Lemma.

Conversely, suppose that $k = k^\prime$ and $d_2 = d_2^\prime + 60\beta$ for some integer $\beta$.  Let $\sigma$ be the filtered ring automorphism on $\bZ \lbrack x \rbrack/(x^4)$ given by
   \[
   \sigma(x) ~=~ x + \beta x^3.
   \]
Since $k = k^\prime$ and $\sigma(x) \equiv x \pmod{x^3}$, it is clear that
   \[
   (\psi^p_S \circ \sigma)(x) ~\equiv~ (\sigma \circ \psi^p_{S^\prime})(x) \pmod{x^3}
   \]
for any prime $p$.  The coefficients of $x^3$ in $(\psi^p_S \circ \sigma)(x)$ and $(\sigma \circ \psi^p_{S^\prime})(x)$ are $(p^2\beta + d_p)$ and $(d_p^\prime + \beta p^6)$, respectively.  They are equal, since
   \[
   d_p - d_p^\prime ~=~ \frac{p^2(p^4 - 1)(d_2 - d_2^\prime)}{60} 
   ~=~ p^2(p^4 - 1)\beta.
   \]
This shows that $\sigma \colon S \to S^\prime$ is an isomorphism of filtered $\lambda$-ring.  This proves the Lemma.
\end{proof}

\medskip
\begin{lemma}
\label{lem3:n=4}
$T$ is isomorphic to one of the sixty $S(k, d_2)$ listed in the statement of Theorem \ref{thm:n=4}.
\end{lemma}

\begin{proof}
The first paragraph of the proof of Lemma \ref{lem1:n=4} also applies to $T$, and so the coefficient of $x^3$ in $\psi^p_T(x)$, call it $d_p$, satisfies \eqref{eq:d_p} for some $k \in \lbrace 1, 3, 5 \rbrace$ and some even integer $d_2$.  Since \eqref{eq:d3} holds in $T$ as well, it follows that $k^2 \equiv 1 \pmod{3}$.  Therefore, $k$ must be either $1$ or $5$.  In other words, $T$ is of the form described in Lemma \ref{lem1:n=4}, and, by Lemma \ref{lem2:n=4}, it is isomorphic to one of the sixty $S(k, d_2)$, as desired.  This proves the Lemma.
\end{proof}

The proof of Theorem \ref{thm:n=4} (2) is complete.
\end{proof}


\subsection{Proof of Theorem \ref{thm:n=4} (3)}
\label{subsec:n=4 (3)}
This case will be proved below as part of Corollary \ref{cor:p>n} (the $n = 4$ case of Conjecture B, which follows from the $n = 3$ case of Conjecture C and Theorem \ref{thm:n>4}.).


\subsection{Proof of Theorem \ref{thm:n=4} (4)}
\begin{proof}
Since the reduction of $R$ modulo $x^3$ is a filtered $\lambda$-ring structure on $\bZ \lbrack x \rbrack/(x^3)$, it follows from Theorem \ref{thm:n=3} (1) that 
   \[
   \psi^p(x) ~=~ \psi^p_R(x) ~\equiv~ 0 \pmod{x^2}
   \]
for all primes $p$.  Thus, $R$ has the form $S((c_p), (d_p))$ as in the statement of Theorem \ref{thm:n=4} (4).

Now let $S = S((c_p), (d_p)) = \lbrace \psi^p(x) = c_p x^2 + d_p x^3 \rbrace$ be a collection of polynomials satisfying $\psi^p(x) \equiv x^p \pmod{p}$.  In order to show that $S$ is a filtered $\lambda$-ring structure on $\bZ \lbrack x \rbrack/(x^4)$, we only need to show that $(\psi^p \circ \psi^q)(x) = (\psi^q \circ \psi^p)(x)$ for all primes $p$ and $q$.  This is true, since both sides are equal to $0$ modulo $x^4$.

To prove the last assertion, let 
   \[
   \bar{S} ~=~ S((\bar{c}_p), (\bar{d}_p)) ~=~ \lbrace \psi^p(x) = \bar{c}_p x^2 + \bar{d}_p x^3 \rbrace
   \]
be another filtered $\lambda$-ring structure on $\bZ \lbrack x \rbrack/(x^4)$.  Suppose first that $S$ and $\bar{S}$ are isomorphic, and let $\sigma \colon S \to \bar{S}$ be an isomorphism.  Write $\sigma(x) = ux + \alpha x^2 + \beta x^3$ with $u \in \lbrace \pm 1 \rbrace$.  Then we have
   \[
   (\psi^p_S \circ \sigma)(x) 
   ~=~ c_p x^2 + (2c_p u\alpha + d_p u)x^3
   \]
and
   \[
   (\sigma \circ \psi^p_{\bar{S}})(x)
   ~=~ u\bar{c}_p x^2 + u\bar{d}_p x^3.
   \]
Conditions (i) and (ii) in the statement of Theorem \ref{thm:n=4} (4) now follow by comparing the coefficients.

The converse is similar.  Indeed, if $(c_p) = u(\bar{c}_p)$ for some $u \in \lbrace \pm 1 \rbrace$ and $\bar{d}_p = d_p + 2c_p \alpha$ for some integer $\alpha$, then the calculation above shows that an isomorphism $\sigma \colon S \to \bar{S}$ is given by
   \[
   \sigma(x) ~=~ ux + \alpha x^2,
   \]
as desired. This proves Theorem \ref{thm:n=4} (4).
\end{proof}

The proof of Theorem \ref{thm:n=4} is complete.


\section{Results about the Conjectures}


\subsection{Proof of Theorem \ref{thm:n>4}}
\begin{proof}
Consider a prime $q$ such that 
   \[
   \vert a_{q,1}^{n-1} - a_{q,1} \vert ~=~ 
   \min\lbrace\vert a_{p,1}^{n-1} - a_{p,1} \vert \colon p \text{ primes}\rbrace.  
   \]
Let $S$ and $T$ be two filtered $\lambda$-ring structures on $\bZ \lbrack x \rbrack/(x^n)$ satisfying \eqref{eq:extension} for all primes $p$, and let $a$ and $b$ be the coefficients of $x^{n-1}$ in $\psi^q_S(x)$ and $\psi^q_T(x)$, respectively.  We claim that if 
   \begin{equation}
   \label{eq1:n>4}
   a ~\equiv~ b \quad \bigl(\text{mod }\vert a_{q,1}^{n-1} - a_{q,1}\vert\bigr),
   \end{equation}
then $S$ and $T$ are isomorphic.  Indeed, if \eqref{eq1:n>4} is true, then we can write
   \begin{equation}
   \label{eq2:n>4}
   a + a_{q,1}c ~=~ b + a_{q,1}^{n-1}c
   \end{equation}
for some integer $c$.  Consider the filtered ring automorphism $\sigma$ on $\bZ \lbrack x \rbrack/(x^n)$ given by 
   \[
   \sigma(x) ~=~ x + cx^{n-1}.
   \]
Then we have that
   \begin{equation}
   \label{eq:pq}
   \begin{split}
   (\psi^q_S \circ \sigma)(x)
   &~=~ \left(\sum_{i=1}^{n-2}\, a_{q,i}(x + cx^{n-1})^i\right) ~+~ a(x + cx^{n-1})^{n-1} \\
   &~=~ \left(\sum_{i=1}^{n-2}\, a_{q,i}x^i\right) ~+~ (a + a_{q,1}c)x^{n-1} \\
   &~=~ \left(\sum_{i=1}^{n-2}\, a_{q,i}x^i\right) ~+~ (b + a_{q,1}^{n-1}c)x^{n-1} \\
   &~=~ (\sigma \circ \psi^q_{T})(x)
   \end{split}
   \end{equation}
Now it follows from \eqref{eq:pq} and \cite[Theorem 4.4]{yau1} that
   \[
   (\psi^p_S \circ \sigma)(x) ~=~ (\sigma \circ \psi^p_T)(x)
   \]
for all primes $p$.  (It is here that we are using the hypothesis $a_{p,1} \not= 0$ for all $p$.)  That is, $\sigma \colon S \to T$ is a filtered $\lambda$-ring isomorphism, as claimed.  It is clear that Theorem \ref{thm:n>4} also follows from this claim.
\end{proof}


\subsection{Proof of Theorem \ref{thm:p>n}}

\begin{proof}
Write $\sigma(x) = \sum_{k=1}^{n-1}\, b_k x^k$ and $\sigma^{-1}(x) = \sum_{i=1}^n\, c_i x^i$.  In particular, we have $b_1 = c_1 \in \lbrace \pm 1 \rbrace$.  Let $p$ be a prime with $p \geq n$.  Write $\psi^p_{\tilde{R}}(x) = \sum_{j=1}^n \, a_j x^j$.  Considering $\psi^p_S(x)$ as a polynomial in $\bZ \lbrack x \rbrack/(x^{n+1})$, we have
   \[
   \begin{split}
   \psi^p_{\tilde{S}}(x) 
   &~=~ (\sigma^{-1} \circ \psi^p_{\tilde{R}} \circ \sigma)(x) \\
   &~=~ \psi^p_S(x) + (b_1^{n+1} a_n + \alpha)x^n,
   \end{split}
   \]
in which $\alpha$ is a sum of integers, each summand having a factor of $a_i$ for some $i = 1, \ldots, n-1$.

If $p > n$, then $\psi^p_{\tilde{R}}(x) \equiv x^p = 0 \pmod{p}$, which implies that $a_i \equiv 0 \pmod{p}$ for $i = 1, \ldots, n$.  Since $\psi^p_S(x) \equiv 0 \pmod{p}$ as well, one infers that $\psi^p_{\tilde{S}}(x) \equiv 0 = x^p \pmod{p}$.

If $p = n$, then we still have $\psi^p_S(x) \equiv 0 \pmod{p}$.  But now $a_p \equiv 1 \pmod{p}$ and $a_i \equiv 0 \pmod{p}$ for $i < p$.  Since $p + 1$ is an even integer, it follows that $\psi^p_{\tilde{S}}(x) \equiv x^p \pmod{p}$, as desired.
\end{proof}


\subsection{Proof of Corollary \ref{cor:p>n}}
\begin{proof}
Consider the case $n = 3$ of Conjecture C first.  In view of Theorem \ref{thm:p>n}, we only need to show that 
   \[
   \psi^2_{\tilde{S}}(x) ~\equiv~ x^2 \pmod{2}.
   \]
Since $\psi^2_S(x) \equiv x^2 \pmod{2}$ already, it suffices to show that the coefficient of $x^3$ in $\psi^2_{\tilde{S}}(x)$ is even.  We will use the same notations as in the proof of Theorem \ref{thm:p>n}.  Observe that the congruence identity, $\psi^2_{\tilde{R}}(x) \equiv x^2 \pmod{2}$, means that $a_1$ and $a_3$ are even and that $a_2$ is odd.  Computing modulo $2$ and $x^4$, we have
   \[
   \begin{split}
   \psi^2_{\tilde{S}}(x) 
   &~\equiv~ \sum_{i=1}^3\, c_i \biggl( \sum_{j=1}^3 \, a_j (b_1x + b_2x^2)^j \biggr)^i \\
   &~\equiv~ \sum_{i=1}^3\, c_i a_2^i(b_1x + b_2x^2)^{2i}.
   \end{split}
   \]
It follows that the coefficient of $x^3$ in $\psi^2_{\tilde{S}}(x)$ is congruent (mod $2$) to $c_1a_2(2b_1b_2) \equiv 0$.  This proves the case $n = 3$ of Conjecture C.

The proofs for the cases $n = 4$ and $5$ of Conjecture C are essentially identical to the previous paragraph.  The only slight deviation is that, when $n = 4$, to check that the coefficient of $x^4$ in $\psi^2_{\tilde{S}}(x)$ is even, one needs to use the fact that $c_2 = -b_1b_2$.

As discussed in the Introduction, Conjectures A and B for $n \leq 6$ follow from Conjecture C for $n \leq 5$ and Theorem \ref{thm:n>4}.
\end{proof}

It should be pointed out that the argument above for the first three cases of Conjecture C does \emph{not} seem to go through for $n = 6$.  In fact, to use the same argument, one would have to check, among other things, that the coefficient of $x^6$ in $\psi^2_{\tilde{S}}(x)$ is even.  But this coefficient turns out to be congruent (mod $2$) to $b_3$.  It may be true that one can always arrange to have an isomorphism $\sigma$ in which $b_3$ is even, but this is only a speculation.  Further works remain to be done to settle the three Conjectures for higher values of $n$.



\begin{thebibliography}{99}
\bibitem{adams1}J. F. Adams, On the non-existence of elements of Hopf invariant one, Ann.\ Math.\ 72 (1960), 20-104.

\bibitem{adams}J. F. Adams, Vector fields on spheres, Ann.\ Math.\ 75 (1962), 603-632.

\bibitem{at}M. F. Atiyah and D. O. Tall, Group representations, $\lambda$-rings and the $J$-homomorphism, Topology 8 (1969), 253-297.

\bibitem{clauwens}F. J. B. J. Clauwens, Commuting polynomials and $\lambda$-ring structures on $\bZ \lbrack x \rbrack$, J. Pure Appl. Algebra 95 (1994) 261-269.

\bibitem{hatcher}A. Hatcher, Algebraic topology, Cambridge Univ.\ Press, 2002.

\bibitem{knutson}D. Knutson, $\lambda $-rings and the representation theory of the symmetric group, Lecture Notes in Math. 308, Springer-Verlag, Berlin-New York, 1973.

\bibitem{mcgibbon}C. A. McGibbon, Self maps of projective spaces, Trans. Amer. Math. Soc. 271 (1982), 325-346.

\bibitem{notbohm}D. Notbohm, Maps between classifying spaces and applications, J. Pure Appl. Algebra 89 (1993), 273-294.

\bibitem{rector}D. L. Rector, Loop structures on the homotopy type of $S^3$, Lecture Notes in Math., vol. 249, Springer-Verlag, Berlin and New York, 1971, 99-105.

\bibitem{wil}C. Wilkerson, Lambda-rings, binomial domains, and vector bundles over $CP(\infty)$, Comm. Algebra 10 (1982), 311-328.

\bibitem{yau1}D. Yau, Moduli space of filtered $\lambda$-ring structures over a filtered ring, Int. J. Math. Math. Sci., (2004), no. 39, 2065-2084.

\bibitem{yau2}D. Yau, On adic genus and lambda-rings, Trans.\ Amer.\ Math.\ Soc., accepted for publication.

\end{thebibliography}
\end{document}